\documentclass[12pt]{article}
 
\usepackage{amssymb,amsmath,amsfonts,amssymb}
\newcommand{\be}{\begin{equation}}
\newcommand{\ee}{\end{equation}}
\newcommand{\la}{\label}
\newcommand{\ba}{\begin{array}{c}}
\newcommand{\ea}{\end{array}}

\renewcommand{\d}{{\mbox{div}}_g}
\newcommand{\D}{ {\Delta_g}}
\newcommand{\nag}{\nabla_g}

\newtheorem{thm}{Theorem}

\title{The Onsager Equation for Corpora}
\author{Peter Constantin\\Department of Mathematics\\The University of Chicago}
\begin{document}
\maketitle 
\begin{abstract} We consider extensions of excluded volume interactions for complex corpora that generalize simple rod-like particles.  The Onsager equation can be defined for quite general configuration spaces, and the dimension reduction of the phase space in the limit of highly intense interaction can be shown.
The formalism describes both freely articulated and interacting N-rods and the example of interacting 2-rods is given in detail.
\end{abstract}
\section{Introduction}
When rod-like particles are suspended in fluids, the particles can be described by relatively simple configuration spaces $M$, for instance, the unit sphere $M= {\mathbb S}^1$ in ${\mathbb R}^2$ or $M={\mathbb S}^2$ in ${\mathbb R}^3$. The probability of finding a rod whose
director $m$ belongs to the region $\Sigma\subset M$ is
$\int_{\Sigma}f(m)dm$. The measure $dm$ on the unit sphere is the uniform measure, the natural volume element (length or area) induced from the ambient Euclidean space. The probability measure $f(m)dm$ characterizes the particle distribution. In equilibrium, this measure is obtained by minimizing the free energy (\cite{ons})
$$
{\mathcal E}[f] = \int_M \left\{\log f + \frac{b}{2}U\right\}fdm
$$
where $b$ is a nonnegative parameter, representing a combination of inverse temperature and  intensity of interaction, and
$$
U(m) = \int_M u(m,p)f(p)dp
$$
is a potential computed from an interaction kernel $u(m,p)$.  The physical modeling is reflected in the choice of kernel. Only two effects are taken into account in the free energy: an entropic effect, representing thermal fluctuation of the particles, and a quadratic mean field effect representing 
excluded volume interactions between particles. 
The description of the problem is completed by specifying the configuration space $M$, the uniform measure $dm$ and the interaction kernel $u$. Once these are determined we may ask questions, such as: do minima of the free energy exist, how many such minima exist, what happens to them as the parameter $b$ varies from $b=0$ to infinity. 

In the case of more complicated corpora, the configuration space $M$ can be 
rather complicated. We use the word ``corpus'' to refer to a body that has finitely many degrees of freedom, such as an assembly of articulated rods, or rods connected to balls of different sizes and with different motion constraints between the parts. The collection of all corpora of a particular problem 
is the configuration space. It is useful to phrase the equilibrium and kinetic problems broadly, in quite general configuration spaces. 
We consider a  metric space $M$, with distance $d$, and assume we are given a Borel probability measure $\mu$ on $M$. This probability measure corresponds to the normalized uniform distribution of corpora $m\in M$, and $M$ is the configuration space.
For absolutely continuous probability Borel measures $\nu<<\mu$, the
potential function
$$
U(x) = \int_M u(x,y)\nu(dy)
$$
associated to $\nu$ is defined by a kernel $u$. The function $u(x,y)$ is given by the model, and it is real valued, symmetric in $x,y$, Lipschitz continuous and bounded. Let $f(x)$ be the Radon-Nikodym density $f = \frac{d\nu}{d\mu}$ so that
$$
\nu = fd\mu,
$$
$f\ge 0, \; \mu-$ a.e., $\int_Mfd\mu = 1$, and let associate to it a free energy
\be
{\mathcal{E}} [f] = \int_M \left\{\log f + \frac{b}{2}U\right\}fd\mu.
\la{freeU}
\ee
Because $\int_MU f d\mu = \int_M\int_M u(x,y)f(x)f(y)d\mu(x)d\mu(y)$, there is no loss of generality by assuming that $u(x,y)$ is symmetric, as only the symmetric part $\frac{1}{2}(u(x,y) + u(y,x))$ contributes to the free energy. Critical points of the free energy obey 
$$
\frac{\delta {\mathcal{E}}}{\delta f} = 0
$$
which leads to Onsager's equation:
\be
f= Z^{-1}e^{-bU}
\la{onsU}
\ee
Indeed, this follows taking the Gateaux derivative (i.e. the first variation)  of ${\mathcal{E}}$ in a direction $h$, i.e.
$$
\frac{d}{dt}{\mathcal{E}}[f+th]_{\left | \right . t=0} = 0,
$$ 
with arbitrary $h$, subject only to $\int_M hd\mu =0$ in order to respect $\int_M(f+th)d\mu =1$. A simple calculation, using the fact that $u(x,y)$ is symmetric in $x,y$ leads to
$\int_M(\log f +b{U}f)hd\mu = 0$ for all $h$. This implies that $\log f +b{U}f$ is constant $\mu$-almost everywhere. Because $u$ is Lipschitz continuous, it is easy to see that the solutions of (\ref{onsU}) are Lipschitz continuous, positive and never vanish.

The kernel $u$ gives the rule to construct the  potential $U$ associated to the particle distribution $\nu = fd\mu$; in examples this is computed from physics and it instructs the particles to conform to each other. As $b\to\infty$ there are phase transitions.  In the abstract case where we have a distance $d$ we may consider general $u(x,y) = \Phi(d(x,y),x,y)$ and 
the simplest examples are $u(x,y) = d(x,y)^p$ with $p=1, 2$ or $u(x,y) = -e^{-\frac{d^2(x,y)}{l^2}}$. The interaction potential 
$u(x,y)$ should have a minimum at $x=y$.

\section{Rods}
The simplest example of single rods has $M={\mathbb S}^1$. The Onsager kernel is $u(\theta_1, \theta_2) = |\sin(\theta_1-\theta_2)|$, and the Maier-Saupe kernel is 
$$
u(\theta_1,\theta_2)  = -\cos^2(\theta_1-\theta_2).
$$ 
Kernels are defined up to addition of constants, so the Maier-Saupe kernel can be taken as $ u(\theta_1, \theta_2) = -\frac{\cos(2(\theta_1-\theta_2))}{2}$.

Throughout this paper we will use the notation 
${\mathcal K}f = -U[f] = -U$ to emphasize the linear dependence of $U$ on $f$. (The minus sign is a convention we maintain to be consistent with previously published literature.)

In the Maier-Saupe potential case, ${\mathcal K}f$ projects $f$ on the eigenspace of the Laplace-Beltrami operator corresponding to eigenvalue 4. If the function $f$ is written as
\be
f(\theta) = \frac{1}{2\pi} + \frac{1}{\pi}\sum_{j=1}^{\infty}y_j\cos(2j\theta)\la{fyj}
\ee
then Maier-Saupe potential is 
\be
{\mathcal K}f = \frac{1}{2}y_1\cos(2\theta).
\la{ms}
\ee
The solutions of the Onsager equation in that case 
are of the form
\be
g(r) = (Z(r))^{-1}e^{r\cos(2\theta)}\la{gr}
\ee
with
\be
Z(r) = \int_0^{2\pi}e^{r\cos(2\theta)}d\theta
\la{zr}
\ee
and in Fourier representation
\be
g(r)(\theta) = \frac{1}{2\pi} + \frac{1}{\pi}\sum_{j=1}^{\infty} g_j(r)\cos(2j\theta)
\la{ggjr}
\ee
Onsager's equation is equivalent to the implicit transcendental equation
\be
g_1(r)) = \frac{2r}{b}.
\la{fixr}
\ee
If $b\le 4$ this equation has one solution, namely $r=0$. If $b>4$ there is
exactly one more solution for $r>0$ at  $r=r(b)$. There is exactly one solution
for $r<0$, at $r = -r(b)$ corresponding to a rotation of $\frac{\pi}{2}$ of the solution with $r>0$. 
It is known also that $g_0(r) =1$ for all $r$,
\be
g_j(r) >0, \;\; g_{j+1}(r) < g_j(r), \quad j\ge 1,\,\, r>0.
\la{pos}
\ee
If $r=r(b)$ is determined by (\ref{fixr}), then 
$$
g_2(r(b)) = 1-\frac{4}{b},
$$
there are recursion relations to compute all $g_j(r)$,
and 
\be
\frac{db}{dr} = \frac{b^2}{2r}\left (g_1^2 - g_2\right )>0
\la{bprime}
\ee
holds for the inverse function, $b(r)$. This is an increasing unbounded function for $b>4$. The stability of the solution $g(r)$ also follows from the inequality $g_1^2-g_2>0$. Both the inequality and the stability are not obvious. The limit of zero temperature is a delta function:
\be
\lim_{b\to\infty} g_{(r(b))}(\theta)d\theta = \delta_0.\la{froze}
\ee
The choice of $\theta=0$ is dictated by the symmetry $\theta\mapsto -\theta$ that we imposed. Otherwise, we obtain any delta function on the circle.
The case of ${\mathbb S}^2$ with Maier-Saupe potential is similar. The potential projects on the eigenspace of eigenvalue $6$ of the Laplace-Beltrami operator on the sphere. Two implicit transcendental equations determine the solution.
There is a finite temperature phase transition and the zero temperature limits are delta functions on points or geodesics. 
\section{Kinetics}
The kinetic description of rod-like particles (\cite{doied}) can be naturally generalized when the configuration space $M$ is a Riemannian manifold.
When $M$ is Riemannian manifold, the kinetic equations are
\be
\partial_t f = \D f - b\d(f\nag({\mathcal K}f))\la{smofp}
\ee
with 
$\D, \d, \nag$  Laplace-Beltrami, divergence and gradient in $M$, ${\mathcal K}f$ given by
\be
{\mathcal K}f = -\int_M u(m,p)f(p)d\mu(p)
\la{k}
\ee
and the uniform measure $d\mu$ is the Riemannian volume element. 
Note that the equation can be written as 
$$
\partial_t f = \d(f\nag(\log f - b{\mathcal K}f))
$$
The solutions are smooth, positive and normalized so they have unit integral. The free energy 
\be
{\mathcal E} = \int \left\{f\log f -\frac{b}{2}f{\mathcal K}f\right\}d\mu
\la{frev}
\ee
is a
Lyapunov functional:
$$
\frac{d}{dt}{\mathcal E} = - \int\limits_{M}f\left| \nag(\log f -b {\mathcal K}f)\right |^2d\mu(p)
$$ 
If $M$ is connected, the only possible steady solutions are solutions of Onsager's equation 
\be
f = Z^{-1}e^{b{\mathcal K}f}.
\la{onsequn}
\ee
The dynamical system is dissipative: the solutions are bounded after an initial transient time. The bounds, in very strong norms, are independent of the initial data. The global attractor is compact, finite dimensional and is formed with solutions of Onsager's equations and their unstable manifolds. 
In the case of $M = {\mathbb S}^1$ with Maier-Saupe potential, the kinetic equations are  given  by the sequence of ODEs
\be
\frac{d}{dt} y_j = -4j^2y_j + bjy_1\left(y_{j-1}-y_{j+1}\right ).
\la{ode}
\ee
The equations on ${\mathbb S}^2$ are more complicated.
Some of the results concerning equilibria and kinetics for the Maier-Saupe potential can be found in \cite{c} -\cite{ckt2}, \cite{ctv}, \cite{cv}, \cite{fatkullin} -\cite{liuzz} and \cite{zar}. Recently it was shown that
the system has inertial manifolds in both ${\mathbb S}^1$ and ${\mathbb S}^2$
(\cite{v1}-\cite{v2}).
\section{Freely articulated corpora}
We consider corpora made of $N$ articulated rods that are allowed to rotate freely. The configuration space is ${\widetilde{M}} = {\mathbb S}^n \times \dots \times {\mathbb S}^n$ with $n =1$ or $2$. The potential is a sum of binary interactions
\be
v(p_1,q_1,p_1,q_2,\dots p_N,q_N) = \sum_{j=1}^N u_j(p_j,q_j).
\la{freev}
\ee
The uniform measure $\mu$ on $\widetilde{M}$ is the product measure. In this situation, the corresponding solution of the Onsager equation is a product measure
\be
\nu(dp_1\dots dp_N) = \Pi_{j=1}^N Z_j^{-1}e^{-bU_j(p_j)}dp_j
\la{pro}
\ee
where each $Z_j^{-1}e^{-bU_j(p_j)}$ is a solution of Onsager's equation in
the $j$ component, and everything can be reduced to the study of individual rods. This situation holds in general, when the corpora are made up of $N$ components, themselves corpora of a simpler kind, each with configuration space $M_j$ and uniform measure $\mu_j$. The composite corpora belong to the configuration space 
$\widetilde{M} = M_1\times \dots \times M_N$, and if the kernel has the form (\ref{freev}) then the corresponding operator
$$
{\widetilde{{\mathcal K}}f}(p_1,\dots p_N) = -\int_{{\widetilde{M}}} v(p_1,q_1,\dots, p_N,q_N)f(q_1,\dots, q_N)d\mu(q_1)\dots d\mu(q_N)
$$
is a sum of corpus-by-corpus operators, corresponding to like-parts only:
$$
{\widetilde {\mathcal{K}}}f(p_1, \dots p_N) = \sum_{j=1}^N \left({\mathcal{K}_j}f\right)(p_j) = -\sum_{j=1}^N U_j[f](p_j)
$$ 
with
$$
{\mathcal{K}_j}f(p_j) = -\int u_j(p_j,q_j)f(q_1,\dots q_N)\Pi d\mu_i(q_i).
$$
In this case, if we consider any solution of
Onsager's equation,  {${\widetilde{f}} = {\widetilde{Z}}^{-1}e^{b{\widetilde{K}}{\widetilde{f}}}$}, then the partition function 
$$
{\widetilde{Z}} = \int_{{\widetilde M}}e^{b{\widetilde K}{\widetilde f}}\Pi d\mu_i
$$
is factored 
$$
{\widetilde{Z}}= \Pi_{j=1}^N Z_j, \quad {\mbox{with}}\; Z_j =\int_{M_j} e^{b{\mathcal K}_jf}d\mu(p_j)
$$
and consequently the solution of Onsager's equation is a product of corpus-by-corpus solutions,
$$
{\widetilde f}(p_1,\dots p_N) = f_1(p_1)f_2(p_2)\dots f_N(p_N)
$$
with $f_j(p_j) = Z_j^{-1}e^{-bU_j(p_j)}$  a solution of the Onsager equation on $M_j$ with uniform measure $\mu_j$ and interaction kernel $u_j$. 
The measure corresponding to the collective is the product measure of the 
single corpus solutions: the components interact only with like-components.
In this case the study of the collective behavior reduces to the study of single corpus statistics, and in particular, the zero temperature limit can be any
combination of zero temperature limits, component-wise.

\section{A simple example: two interacting rods}
Clearly, if the components of complex corpora do not interact only with like-components, then a product measure is not the solution of the Onsager equation. Let us start by investigating a simple example. Assume that two rods with configuration space $M= {\mathbb{S}^1}$ are articulated, so together the corpora belong to
${\widetilde{M}} = {\mathbb S}^1\times {\mathbb S}^1$. We will assume that the interaction between pairs of two-rods is determined entirely by the area of the triangle formed by each two-rod, and each two-rod corpus rejects two-rods that have very different area than themselves. The simplest interaction that achieves this is
$$
u(p_1,q_1, p_2, q_2) = \|e(p_1)\wedge e(p_2) -e(q_1)\wedge e(q_2)\|^2
$$ 
with $e(p) = (\cos p, \sin p)$ if $p \in [0, 2\pi]$. Each ${\mathbb{S}}^1$
is viewed as a subset of ${\mathbb R}^2$; the exterior product $\bigwedge^2 {\mathbb R}^2$ is isomorphic to ${\mathbb R}$. Every element in it is a multiple of $e_1\wedge e_2$, with $e_1 = (1,0)$ and $e_2 = (0,1)$. For instance, $e(p_1)\wedge e(p_2) = \sin(p_2-p_1)(e_1\wedge e_2)$. The norm satisfies
$$
\|e(p_1)\wedge e(p_2) - e(q_1)\wedge e(q_2)\|^2 = \left (\sin(p_1-p_2) - \sin(q_
1-q_2)\right)^2
$$
and consequently, the integral operator is
$$
{\mathcal K}f(p_1, p_2) = -\int_{{\widetilde {M}}} \|e(p_1)\wedge e(p_2) - e(q_1)
\wedge e(q_2)\|^2 f(q_1,q_2)dq_1dq_2.
$$ 
We take the uniform measure on ${\mathbb S}^1\times{\mathbb S}^1$ to be the natural measure $dq_1 dq_2$.
We note that the potential has the form 
$$
U(p_1, p_2) =  \sin^2(p_1-p_2) -2z\sin(p_1-p_2) + \gamma
$$
with $z, \gamma$ constants determined by $f$. Onsager's equation $f = Z^{-1}e^{b{\mathcal K f}}$ reduces therefore to
$$
\left\{
\ba
z = [\sin\theta](z,\gamma)\\
\gamma = [\sin^2\theta](z,\gamma)
\ea
\right .
$$
where we use the notation 
$$
[\phi](z,\gamma) = \int_0^{2\pi}\phi(\theta)g(\theta)d\theta
$$
$$
g(\theta) = Z^{-1}e^{-b\sin^2(\theta) + 2bz\sin\theta - b\gamma}
$$
$$
Z= \int_0^{2\pi}e^{-b\sin^2(\theta) +2 bz\sin\theta - b\gamma}d\theta
$$
The solution is $f(p_1, p_2) = g(p_1-p_2)$. Note that $g$ does not depend on $\gamma$. Let then 
$$
s(\theta, z) = \sin\theta -z,
$$
and let
$$
[s](b,z) =\frac{\int_0^{2\pi}s(\theta,z)e^{-bs^2(\theta,z)}d\theta}{\int_0^{2\pi}e^{-bs^2(\theta, z)}d\theta}.
$$
The Onsager equation is equivalent to
\be
[s](b,z) =0.
\la{onsu}
\ee
This determines $z$, which in turn determines $g$, $f$. Note that $z=0$ always a solution that yields
$$
f_0(p_1,p_2) = Z^{-1}e^{-b\sin^2(p_1-p_2)}.
$$
As $b\to\infty$ this solution tends to $\delta ((p_1-p_2)\,\mbox{mod}\pi)$, a degenerated two-rod.

Consider now 
$$
\lambda(z,\tau) = b^{\frac{1}{2}}\int_0^{2\pi}e^{-b(\sin\theta -z)^2}d\theta
$$
with  $\tau = b^{-1}$.
Note that
$$
[s](\tau^{-1}, z) = \frac{1}{2b}\frac{\partial_z\lambda}{\lambda}
$$
so (\ref{onsu}) is equivalent to
\be
\partial_z \lambda  = 0.
\la{onsl}
\ee
Note also that $\lambda$ solves a linear heat equation with temperature 
as time:
$$
\partial_{\tau}\lambda = \frac{1}{4}\partial_z^2\lambda
$$
The function $\lambda$ is even in $z$, so it is enough to study it on $[0,1]$. The ``initial value'' obeys, for any $z\in [0,1)$,
$$
\lim_{\tau\to 0}\lambda(z,\tau) = 2{\sqrt{\pi}}\frac{1}{\sqrt{1-z^2}}
$$
This is an increasing function for $z\in (0,1)$, the derivative is positive. Because the derivative obeys the same linear heat equation, and
the heat equation preserves positivity, we might have been tempted to think that $z=0$ is the only solution of $\partial_z\lambda = 0$. But the behavior is more subtle because, clearly, by direct inspection, $\frac{\partial\lambda}{\partial z}(1, \tau)<0$! The reason for this is the singular behavior at $z=1$. As $\tau\to 0$ or $b\to\infty$, we have a change of behavior of the function in a small region of the order $b^{-1/2}$ near $z=1$. In this region there is a transition to much higher values of $\lambda$, and the derivative $\partial_z \lambda$ changes sign. This implies that there is a phase transition at positive $\tau$, i.e.,  there exists $0<z(b)$ satisfying
$$
\partial_z\lambda (z_b,\tau) = 0
$$
and $\lim_{\tau\to 0}z_b=1$. Consequently 
$$
\lim_{b\to\infty}f(p_1-p_2) =\delta\left (\left (p_1-p_2 -\frac{\pi}{2}\right ){\mbox{mod}}\pi\right)
$$
is a delta function concentrated on right two-rods that make a right angle.
It is instructive to note that the reason for this transition is that the
vanishing of the gradient of the phase $\sin\theta -z$ can occur, when
$z=1$ at a minimum of the potential. Further details will presented elsewhere (\cite{cz}).

\section{A few general observations}
More complicated configuration spaces arise when the corpora are  n-gons in space or in the plane. The natural conformation
distance between such corpora is Hausdorff distance, modulo rotations. This area is rather open for investigations, and it is useful to start by stating the main general expectations. It is expected that, if the configuration spaces are compact (a realistic assumption if the corpora have finitely many degrees
of freedom and finite extensivity), then generically, the zero temperature limit will be a singular measure, i.e., it will be concentrated on a set of zero
$\mu$ measure, where $\mu$ is the uniform measure. In many examples the set of zero measure is just one point in $M$.
 
Let us consider a general  compact metric space $M$ with distance $d$ 
and ``uniform'' distribution $\mu$, a Borel probability on $M$. 
Let us assume that there exist $0<k<1$, $c>0$
\be
\mu(B(x,r))\ge c e^{-r^{-k}}
\la{asu}
\ee
for all $x\in M$, and all $r$ sufficiently small. Here $B(x,r)$ is the ball centered at $x\in M$ and of radius $r$ in the metric $d$. This assumption says that all balls are charged at least a small positive amount.  Note that if $M$ is any compact Riemannian manifold of dimension $n$ and $\mu$ is the Riemannian volume element, then the condition 
(\ref{asu}) is automatically satisfied, because $r^n$ is much larger, for small $r$, than the right hand side of (\ref{asu}).
The interaction kernel we consider is a function of the distance, $u(d(x,y))$, so
\be
U(x) = \int_M u(d(x,y))f(y)d\mu(y).
\la{Uf}
\ee
We assume that $u$ is non-negative, bounded and Lipschitz continuous,
i.e. there exist positive constants $C$ and $L$ so that
\be
0\le u(d) \le C
\la{ub}
\ee
and 
\be
\left |u(d_1) -u(d_2)\right | \le L\left |d_1-d_2\right|
\la{lipu}
\ee
hold for all $d, d_1, d_2 \ge 0$. We assume also that
\be
u(0) =0.
\la{uzero}
\ee
As we mentioned before, the interaction kernels are defined up to additive constants: if we add 
$c$ to the interaction kernel, then the potential is changed  by the same amount $c$ and the free energy is changed by adding $\frac{bc}{2}$; its critical points, and in particular its minima are unchanged.
Let us take a ball of $B$ of radius $r$, set $\chi = \mu(B)^{-1}{\mathbf 1}_B$ the normalized indicator function of $B$, and compute
$$
{\mathcal E}[\chi] = \log(\mu(B)^{-1}) + \frac{b}{2}\mu(B)^{-2}\int_B\int_Bu(d(x,y))d\mu(x)d\mu(y)
$$
We obtain, using (\ref{asu}, \ref{ub}, \ref{lipu}) that
$$
{\mathcal E}[{\chi}] \le \left(\frac{1}{r}\right)^k + \log\left(\frac{1}{c}\right) + bLr
$$
which implies, by choosing $r=b^{-1}$, that
$$
\inf_f {\mathcal E}[f] \le b^k + C_1
$$
holds for large $b$ with some constant $C_1$. This means that uniform measure, whose energy is linear in $b$,
$$
{\mathcal E}[1] = \frac{b}{2}\int_M\int_M u(d(x,y))d\mu(x)d\mu(y)
$$
does not achieve the minimum of energy for large $b$. On the other hand, for may examples, the uniform measure is a solution of Onsager's equation. Indeed,
if $T$ is a $\mu$ measure preserving isometry of $M$, then the potential associated to the uniform measure $\mu$, $U[1] = -{\mathcal K}1$, is invariant under right composition with $T$, i.e.,
$$
U(Tx) = U(x)  
$$
holds for all $x\in M$ and all measure preserving isometries. If measure preserving isometries act transitively, i.e., for any $x,y\in M$ there exists $T$
a measure preserving isometry that maps $x$ to $y$, $y=Tx$ then, $U = -{\mathcal K}(1)$ is a constant. This implies that the uniform measure is a solution of
Onsager's equation. This is the case for many homogeneous spaces.
The combination of these two very simple observations leads to the conclusion that a phase transition occurs, whenever the homogeneous measure is a solution of Onsager's equation, the condition (\ref{asu}) holds and the interaction kernel is a normalized Lipschitz function of distance. That simply means that, while for $b=0$ obviously the only solution of Onsager's equation is the homogeneous measure, and while this continues to be a solution for $b>0$, at large enough $b$ there exist other solutions as well, in quite great generality.

Now we describe a general tendency of solutions to 
concentrate. 
Let $d\nu = fd\mu$ be any probability measure absolutely continuous with respect to $\mu$ and let $U$ be the potential associated to it via (\ref{Uf}). Then, in view of the property (\ref{ub}) we have
\be
0\le U(x) \le C\la{Ub}
\ee
and from (\ref{lipu}) and the triangle inequality we deduce that
\be
\left |U(x) -U(y)\right | \le Ld(x,y)
\la{lipU}
\ee
so that the potentials associated to any probability $fd\mu$ are non-negative, uniformly bounded and Lipschitz continuous.
\begin{thm} Let $M$ be a compact metric space with distance $d$. Let $\mu$ be a Borel probability measure on $M$ that satisfies (\ref{asu}). Let $u$ satisfy
(\ref{ub},\ref{lipu}). Then:

(A) For any $b>0$ there exists a solution $g$ that minimizes the enery:
$$
{\mathcal E}[g] = \min_{f>0,\; \int_Mfd\mu =1}{\mathcal E}[f]
$$
The function $g$ solves the Onsager equation
$$
g(x) = (Z(b))^{-1}e^{-bU(x)}
$$
with 
$$
Z(b) = \int_Me^{-bU(x)}d\mu(x)
$$
and
$$
U(x) = \int_M u(d(x,y))g(y)d\mu(y)
$$
The function $g$ is normalized $\int gd\mu =1$, strictly positive and Lipschitz continuous.

(B) Let $b_n\to\infty$ and let $d\nu_n = g_nd\mu$ be a sequence of solutions of
Onsager equations corresponding to $b_n$. By passing to a subsequence we may assume that the sequence converges weakly to a probability measure $\nu= \lim_n\nu_n$. There exists a non-negative Lipschitz continuous function $U_{\infty}(x)$ on $M$ such that $\nu$ is concentrated on the set
$$
\Sigma = \{x\in M\; \left | \right. U_{\infty}(x) = {\mbox{min}}_{y\in M}U_{\infty}(y)\}
$$
Thus, for any continuous function $\phi$ supported in the open set $M\setminus \Sigma$,
$$
\lim_{n\to\infty}\int_M\phi(x)g_n(x)d\mu = 0
$$
\end{thm}
\noindent{\bf Proof.} For the proof of (A) we fix $b>0$ and 
note that ${\mathcal E}[f]$ is bounded below uniformly for all
$f>0, \; \int_Mfd\mu =1$. We take then 
a minimizing sequence $f_j$, 
$$
a = \inf_{f>0; \int_M fd\mu =1}{\mathcal E}[f] = \lim_{j\to\infty}{\mathcal E}[f_j].
$$
Without loss of generality, by passing to a subsequence and relabelling, we may assume that the measures $f_jd\mu$ converge weakly to a measure $d\nu$.
Using (\ref{Ub}, \ref{lipU}) and the Arzela-Ascoli theorem, we may pass to a subsequence, which we relabel again 
$f_{j}$, so that $U_{j}$ converge uniformly to a non-negative Lipschitz continuous function $U$. Then it follows that 
$$
U(x) = \int_M u(d(x,y))d\nu(y)
$$
holds and 
$$
\lim_{j\to\infty}\int_M U_j f_j d\mu = \int_M Ud\nu.
$$
Because ${\mathcal E}[f_j]$ is a convergent sequence of numbers, it follows that
$$
\lim_{j\to\infty}\int_Mf_j\log f_j d\mu
$$
exists. In particular, the above integrals are bounded uniformly. 
It then follows that $d\nu$ is absolutely continuous i.e., $d\nu = gd\mu$,  
with $g\ge 0$, $g\in L^1(d\mu)$. Indeed, the sequence $f_jd\mu$ is uniformly absolutely continuous. This follows from the convexity of the function $y\log y$ and the Jensen inequality
$$
(\mu(A))^{-1}\int_A f\log f d\mu \ge m\log m
$$
where $m = \mu(A)^{-1}\int_A fd\mu$. Then we have
$$
m\log m\le C/\mu(A)
$$
with a fixed constant $C>0$, uniformly for all $f= f_j$ and any $A$.
Let us choose $R$ so that $R\log R = C/\mu(A)$. If we denote by $I = \int_Afd\mu$, then either $m\le R$, or, if not, then $m\log R\le C/\mu(A)$. In either case the inequalities imply
$$
I \le C/\log (R)
$$  
and as $\mu(A)\to 0$, $R\to\infty$. This inequality signifies
$$
\int_A f_nd\mu \le \delta(\mu(A))
$$
with $\lim_{x\to 0}\delta(x) = 0$ and $\delta(x)$ independent of $n$. This implies that $\nu$ is absolutely continuous, $\nu = gd\mu$ with $0\le g\in L^1(d\mu)$. The weak convergence tested on the function $1$ implies that $\int gd\mu =1$. In general, weak convergence of measures is not enough to show lower semicontinuity of nonlinear integrals or almost everywhere convergence. We claim however that in fact the convergence $f_n\to g$ takes place strongly in $L^1(d\mu)$:
$$
\lim_{n\to\infty}\int_M |f_n(x) - g(x)|d\mu(x) = 0.
$$
In order to prove this we prove that $f_n$ is a Cauchy sequence in $L^1(d\mu)$.
We take $\epsilon>0$ and choose $N$ large enough so that
$$
\sup_{x\in M}\left |U_n(x) -U(x)\right| \le \frac{\epsilon^2}{10 b},
$$
$$
\left |\int_M U(x)(f_n(x)- f_m(x))d\mu\right| \le \frac{\epsilon^2}{10 b}
$$
and
$$
{\mathcal E}[f_n] \le a + \frac{\epsilon^2}{16}
$$ 
hold for $n,m\ge N$. Let $s(x) = \frac{1}{2}(f_n(x) + f_m(x))$. Then $\int_M sd\mu =1$,  $s>0$, so
$$
a\le {\mathcal E}[s].
$$
Therefore
$$
\frac{1}{2}\left \{{\mathcal E}(f_n) + {\mathcal E}(f_m)\right\} - {\mathcal E}[s] \le \frac{\epsilon^2}{16}.
$$ 
On the other hand,
$$
\int_M \left \{\frac{1}{2}\left [f_n\log f_n + f_m\log f_m\right ] - s\log s\right\}d\mu \le 
\frac{1}{2}\left \{{\mathcal E}(f_n) + {\mathcal E}(f_m)\right\} - {\mathcal E}[s] + \frac{\epsilon^2}{16}
$$
so
$$
\int_M \left \{\frac{1}{2}\left [f_n\log f_n + f_m\log f_m\right ] - s\log s\right\}d\mu \le \frac{\epsilon^2}{8}.
$$
Denote $\chi = \frac{f_n-f_m}{f_n + f_m}$ and note that $-1\le \chi\le 1$ holds $\mu$ - a.e. Also, elementary calculation show that
$$
\left \{\frac{1}{2}\left [f_n\log f_n + f_m\log f_m\right ] - s\log s\right\} =
\frac{s}{2} G(\chi)
$$
holds with
$$
G(\chi) = \log(1-\chi^2) + \chi\log\left(\frac{1+\chi}{1-\chi}\right).
$$
Note that $G$ is even on $(-1,1)$, that $G^{\prime}(\chi) = \log\left(\frac{ 1+\chi}{1-\chi}\right)$, $G(0) = G^{\prime}(0) = 0$ and $G^{\prime\prime}(\chi) = \frac{2}{1-\chi^2}\ge 2$ on $(-1,1)$. Consequently,
$$
0\le \chi^2 \le G(\chi)
$$
holds for $-1\le \chi \le 1$. It follows that we have 
$$
\int_M \frac{(f_n-f_m)^2}{f_n+f_m} d\mu \le \frac{\epsilon^2}{2}
$$
But, writing $|f_n-f_m| = \sqrt{f_n+f_m}\frac{|f_n-f_m|}{\sqrt{f_n+f_m}}$
and using the Schwartz inequality we deduce
$$
\int_M|f_n-f_m|d\mu \le \epsilon.
$$
Therefore the sequence $f_n$ is Cauchy in $L^1(d\mu)$. This proves that
the weak limit $f_nd\mu \to gd\mu$ is actually strong $f_n\to g$ in $L^1(d\mu)$. By passing to a subsequence if necessary, we may assume also that $f_n\to g$ holds also $\mu$- a.e. Then, from Fatou's Lemma
$$
\int_Mg\log gd\mu \le \lim_{j\to\infty}\int f_j\log f_j d\mu.
$$
This implies that
$$
{\mathcal{E}}[g] = a.
$$
The fact that $g$ solves the Onsager equation follows by taking the Gateaux derivative, and thus
$$
g = Z^{-1}e^{-bU}
$$
with $Z = \int_M e^{-bU}d\mu$. Because $U$ is bounded it follows that
$g$ never vanishes and because $U$ is Lipschitz, so is $g$.

The proof of (B). Let $b_n\to\infty$ and let us take a subsequence so that
$g_nd\mu$ converges weakly to $d\nu$. As above, because of  (\ref{Ub}, \ref{lipU}) and the Arzela-Ascoli theorem, we may pass to a subsequence, which we relabel again
$g_{n}$, so that $U_{n}$ converge uniformly to a non-negative Lipschitz continuous function $U_{\infty}$. Let $x_n$ be a point where $U_n(x)$ attains its minimum $U_n(x_n) = \min_{x\in M}U_n(x)$. By passing again to a subsequence we may assume that $x_n$ converge to some point $x$. It follows that that $U_{\infty}(x) = \min_{m\in M}U_{\infty}(m) =\alpha$. Let $\phi$ be a continuous function in $M$ compactly supported in $M\setminus \Sigma$ where $\Sigma = \{m\in M\;\left |\right . U_{\infty}(m) = \alpha\}$. There exists $\epsilon>0$ so that, for every $m$ in the support of $\phi$, $U_{\infty}(m)\ge \alpha + 4\epsilon$. Let us take $N$ so large that
$\sup_M\left |U_{\infty}(m)-U_n(m)\right | \le \epsilon$ for $n\ge N$ and $d(x_n,x)\le \frac{\epsilon}{L}$ Denote $\alpha_n$the minimum of $U_n$. It follows that 
$|\alpha_n-\alpha|\le 2\epsilon$ and $U_n(m) \ge \alpha_n + \epsilon$ on the support of $\phi$.
On the other hand, we have    
$$
Z_n(b_n) \ge \int_{B(x_n, \frac{1}{b_n})}e^{-b_nU_n(z)}d\mu(z)\ge e^{-b_n\alpha_n}\int_{B(x_n,\frac{1}{b_n})}e^{-b_nLd(z,x_n)}d\mu(z),
$$
and using (\ref{asu}) we get
\be
Z_n(b_n)\ge e^{-L}ce^{-b_n^{k}}e^{-b_n\alpha_n}
\la{lz}
\ee
Therefore, on the support of $\phi$ we have
$$
g_n(m) = (Z_n(b))^{-1}e^{-b_n U_n(m)}\le e^L c^{-1}e^{b_n^k}e^{b_n\alpha_n}e^{-b_n(\alpha_n +\epsilon)}
$$
and consequently
\be
\left |\int_{M}\phi(m)g_n(m)d\mu(m)\right | 
\le e^L{c}^{-1}e^{b_n^k}e^{-b_n\epsilon}\int_M|\phi(m)|d\mu(m)
\la{conc}
\ee
holds, and therefore, as $k<1$ we have
$$
\lim_{n\to\infty}\int_M\phi g_nd\mu = 0.
$$
\section{Conclusions} The generalization of excluded volume interactions, from simple rod-like particles to complicated corpora leads to a ``simple'' equation
(the Onsager equation) in ``complicated'' spaces. The examples of single rods, and articulated two-rods show significant complexity reduction. The complexity reduction, once the problem is phrased correctly, is expected to be generic.
The zero temperature or high intensity limit of probability distributions of 
corpora  concentrates on the minima of certain Lipschitz functions, in general.   

\section{Acknowledgments} Work partially supported by NSF-DMS grant 0504213.
I thank A. Zlatos and B. Farb for fruitful conversations.  
\bibliographystyle{amsalpha}

\end{document}